\documentclass[12pt]{article}

\title{Computing expectiles via fixed point iterations}

\author{
Ha Thi Khanh Linh\footnote{Free University of Bolzano-Bozen, \href{mailto:thikhanhlinh.ha@unibz.it}{thikhanhlinh.ha@unibz.it}} \;
Andreas H. Hamel\footnote{Free University of Bolzano-Bozen, \href{mailto:andreas.hamel@unibz.it}{andreas.hamel@unibz.it}}\;
Daniel Kostner\footnote{formerly Free University of Bolzano-Bozen, \href{mailto:daniel.kostner@mailbox.org}{daniel.kostner@mailbox.org}}
}
\date{{\small \today}}

\usepackage{array} 
\usepackage{verbatim} 

\usepackage{amsmath, amssymb, enumerate, color}
\usepackage{bbm}
\usepackage{graphicx,float}

\newtheorem{theorem}{Theorem}
\newtheorem{corollary}[theorem]{Corollary}
\newtheorem{remark}[theorem]{Remark}

\newtheorem{proposition}[theorem]{Proposition}
\newtheorem{example}[theorem]{Example}

\numberwithin{equation}{section}  
\numberwithin{figure}{section}    
\numberwithin{table}{section}     
\numberwithin{theorem}{section}

\newcommand{\norm}[1]{\ensuremath{\left\| #1 \right\|}}
\newcommand{\abs}[1]{\ensuremath{\left| #1 \right|}}
\newcommand{\cb}[1]{\ensuremath{ \left\{ #1 \right\} }}

\newcommand{\bs}{\backslash}

\newcommand{\pend}{ \hfill $\square$ \smallskip}

\newcommand{\R}{\mathrm{I\negthinspace R}}
\newcommand{\N}{\mathrm{I\negthinspace N}}
\newcommand{\E}{\mathbb{E}}

\DeclareMathOperator*{\argmin}{argmin}

\newcommand{\One}{\mathrm{1\negthickspace I}}

\usepackage[
colorlinks=true, linkcolor=blue, citecolor=blue, urlcolor=blue, filecolor=blue
auskommentieren und unten pdfborder aktivieren]{hyperref}

\definecolor{color0}{gray}{.50}
\definecolor{color1}{rgb}{0,.2,.8}
\definecolor{color2}{rgb}{1,.2,0}
\definecolor{color3}{rgb}{.2,.7,.6}

\begin{document}

\maketitle

\begin{abstract}
Expectiles are statistical parameters which also provide a class of sublinear risk measures in finance. They are solutions of continuous optimization problems. The corresponding first order condition provides two different fixed point characterizations for expectiles, both of which can be utilized for computing them. Although especially the so-called two-sided version is already implemented and widely used, a general convergence proof appears to be new.
\end{abstract}



\section{Introduction}

Expectiles as statistical parameters, basically M-quantiles, were introduced within a regression framework by Newey and Powell in \cite{NeweyPowell87Economet}. Since then, they have gained considerable attention supported by the fact that they can be turned into financial risk measures which are coherent, law-invariant and elicitable \cite{BelliniDiBernardino17EJF, BelliniEtAl14IME, Ziegel16MF}, and they are the only ones with these features within the class of M-quantiles.

The computation of expectiles is a non-trivial task. We refer the reader to the recent paper by Daouia, Stupfler and Usseglio-Carleve \cite{DaouiaStupflerUsseglio24SC} which not only contains some exact computation formulas for the discrete case and numerical procedures such as a Newton type algorithm but also many examples illustrating several features and difficulties. Results on the asymptotic behavior of expectiles, i.e., the convergence of sample expectiles to the expectile of the underlying distribution can be found in \cite{HolzmanKlar16EJS}.

In the present note, the following results are given:

-- a direct proof for the convergence of a one-sided fixed point iteration via Banach's fixed point theorem, 

-- a direct proof of a two-sided fixed point iteration under less restrictive assumptions than the ones in \cite[Section 2.2]{DaouiaStupflerUsseglio24SC},

-- ramifications for the discrete case including finite termination results.

The first and the second contribution seem to be new in the literature; in particular, we could not find a general convergence proof for the fixed point iteration used in \cite{DaouiaStupflerUsseglio24SC} and also already in the software package \cite{SobotkaEtAl14Github}. One may note that the computation of univariate expectiles is also important for the computation of projection expectiles in the multivariate case without \cite{CascosOchoa21JMVA} or with underlying order relation for the values of random vectors \cite{HaHamel25JMVA}. In \cite{BelliniDiBernardino17EJF}, the method for the computation of expectiles was not indicated.

\section{Computation of univariate $\alpha$- expectiles}

\subsection{Fixed point iterations}

The $\alpha$-expectile $e_{\alpha}(X)$ of a random variable $X \colon \Omega \to \R$ with existing second moment for $\alpha \in (0, 1)$ is the unique minimizer of an asymmetric quadratic loss function:
\begin{equation}
\label{EqMinProblem}
e_{\alpha}(X) = \argmin_{x \in \R} \left \{(1-\alpha)\E(X-x)_{-}^{2} + \alpha \E(X-x)_{+}^{2}\right \}.
\end{equation}
The first order condition for the minimization problem yields the equation
\begin{equation} 
\label{EqFOC}
\alpha \E(X - x)_+ = (1-\alpha)\E(X - x)_-
\end{equation}
which can also be used to extend the definition of the $\alpha$-expectile to random variables with (only) first moments: it is the unique number $x = e_{\alpha}(X)$ satisfying \eqref{EqFOC}.

Replacing the positive and the negative part in \eqref{EqFOC} by the other one, respectively, leads to the ''one-sided'' formulas (see \cite[formula (12)]{BelliniEtAl14IME},  \cite[formula (3)]{CascosOchoa21JMVA} and already \cite{NeweyPowell87Economet}) 
\begin{align}
\label{EQOneSidedFP-A}
x & = \E X + \frac{2\alpha - 1}{1-\alpha} \E(X - x)_+,  \\[.2cm]
\label{EQOneSidedFP-B}
x &  = \E X + \frac{2\alpha - 1}{\alpha} E(X - x)_-
\end{align}
which could be used for a fixed point iteration. The following result shows that this idea is justified. 

\begin{theorem}
\label{ThmContraction} The function $\varphi_\alpha \colon \R \to \R$ defined by
\[
\varphi_\alpha(x) =  
\left \{ 
	\begin{aligned} 
	& \E X + \frac{2\alpha - 1}{1-\alpha} \E(X - x)_+ \quad \text{for} \quad 0 < \alpha < \frac{1}{2}\\[.2cm]
	& \E X + \frac{2\alpha - 1}{\alpha} \E(X - x)_- \quad \text{for} \quad \frac{1}{2} < \alpha < 1 
	\end{aligned} 
\right.
\]
is a contraction.
\end{theorem}

{\sc Proof.} The functions $x \mapsto (y-x)_+$, $x \mapsto (y-x)_-$ are Lipschitz continuous with constant 1 for all $y \in \R$. This property is inherited by $x \mapsto \E(X - x)_+$ and $x \mapsto \E(X - x)_-$. One has $0 < \abs{\frac{2\alpha - 1}{1-\alpha}} = \frac{1- 2\alpha}{1-\alpha} < 1$ for $0 < \alpha < \frac{1}{2}$, thus $\varphi_\alpha$ is a contraction with contraction constant $\frac{2\alpha - 1}{1-\alpha}$ if $0 < \alpha < \frac{1}{2}$. Similarly, $\frac{1}{2} < \alpha \leq 1$ implies $0 < \abs{\frac{1 - 2\alpha}{\alpha}} =  \frac{2\alpha - 1}{\alpha} < 1$ which is the contraction constant in this case. \pend

Banach's fixed point theorem \cite[3.48 Contraction Mapping Theorem]{AliprantisBorder06Book} now yields the convergence of a sequence defined by $x_{i+1} = \varphi_\alpha(x_i)$ with any $x_0 \in \R$ to the unique fixed point of $\varphi_\alpha$, i.e., $e_\alpha(X)$. Moreover, $\varphi_{0.5}(x) = \E X = e_{0.5}(X)$.

The fixed point version of \eqref{EqFOC} is also mentioned in \cite{DaouiaStupflerUsseglio24SC}, namely for $\varphi_\alpha$ with $1/2 < \alpha < 1$, but the fixed point iteration is neither mentioned, nor applied in this reference.
 
\begin{remark}
\label{RemFixedPointElicitability}
Theorem \ref{ThmContraction} can also be used to give another definition of the $\alpha$-expectile: it is the unique fixed point of the function $\varphi_\alpha$. This is a different characterization compared to elicitability \cite{Ziegel16MF} which raises questions like the following. Is there a useful fixed point characterization of quantiles? Can a fixed point characterization of such a statistical parameter be used for purposes like backtesting? We note, however, that it might be difficult to turn \cite[Corollary 3, Equation (6)]{BelliniEtAl14IME} into a fixed point form without further assumptions to the functions $\Phi_1$, $\Phi_2$ defining the objective of the minimization problem parallel to \eqref{EqMinProblem}.
\end{remark}

Using
\begin{equation}
\label{EqPosNegDiff}
(X - x)_+ = (X - x)\One_{X>x}, \quad (X - x)_- = (x - X)\One_{X \leq x}
\end{equation}
one can  give \eqref{EqFOC} the ''two-sided'' form
\begin{equation}
\label{EqTwoSidedFP}
x = \psi_\alpha(x) := \frac{\alpha \E X\One_{X > x} + (1-\alpha)\E X\One_{X \leq x}}{\alpha \E\One_{X > x} + (1-\alpha)\E\One_{X \leq x}}.
\end{equation}
which in turn means that the expectile is a fixed point of the function $\psi_\alpha(x)$.

 A fixed point iteration for $\psi_\alpha$ is the underlying idea of the Repeated Weighted Averaging (RWA) algorithm used in \cite{SobotkaEtAl14Github} for computing sample expectiles. In fact, the recursive formula in \cite[p. 7]{DaouiaStupflerUsseglio24SC}     
 \[
 x_{i+1} = \frac{(2\alpha - 1)\mathbb{E}(X \mathbbm{1}_{\{X > x_i\}}) + (1 - \alpha)\mathbb{E}(X)}{(2\alpha - 1)\mathbb{P}(X > x_i) + 1 - \alpha}
 \]
 (a version of the fixed point iteration for $\psi_\alpha$) is also used by the codes in \cite{SobotkaEtAl14Github} to compute the empirical expectiles. This was probably not realized in \cite{DaouiaStupflerUsseglio24SC} since it was not stated in the description/manual of \cite{SobotkaEtAl14Github}. On the other hand, \cite{DaouiaStupflerUsseglio24SC} gives an interpretation of the fixed point iteration for $\psi_\alpha$ as a Newton-Raphson method under additional assumptions.

To the best of our knowledge, there is no general convergence proof yet for the fixed point iteration via \eqref{EqTwoSidedFP} although the algorithm works very well. Under additional assumptions (the distribution has a continuous density function), \cite[Thm. 2.2, Cor. 2.3]{DaouiaStupflerUsseglio24SC} provides quadratic convergence as expected from a Newton-type method. 

The remainder of this section is devoted to a convergence proof for the sequence defined by $x_{i+1} = \psi_\alpha(x_i)$ with $x_0 \in \R$ without any further assumption. Define the functions $\gamma_\alpha, h_\alpha \colon \R \to \R$ by
\[
h_\alpha(x) = \alpha \E(X - x)_+ - (1-\alpha)\E(X - x)_-
\]
and
\[
\gamma_\alpha(x) = \alpha\overline F_X(x) + (1-\alpha) F_X(x).
\]
observing $\E\One_{X > x} = \overline F_X(x)$, $\E\One_{X \leq x} = F_X(x)$ where $F_X$ is the cdf of $X$ and $\overline F_X = 1- F_X$ its survival function.
Of course, $0 \leq \overline F_X(x) \leq 1$ for all $x \in \R$ and 
\[
\min\{\alpha, 1-\alpha\} \leq \gamma_\alpha(x) \leq \max\{\alpha, 1-\alpha\}
\]
due to the fact that $\gamma_\alpha(x)$ is a convex combination of $\alpha$ and $1-\alpha$.

\begin{proposition}
\label{PropPsiRep}
(a) One has
\begin{equation}
\label{EqTwoSidedFP-Mod}
\psi_\alpha(x) - x = \frac{h_\alpha(x)}{\gamma_\alpha(x)} = \frac{\alpha \E(X - x)_+ - (1-\alpha)\E(X - x)_-}{\alpha \E\One_{X > x} + (1-\alpha)\E\One_{X \leq x}}.
\end{equation}

(b) The function $h_\alpha$ is decreasing.

(c) If $x \leq e_\alpha(X)$, then $0 \leq h_\alpha(x)$; if $e_\alpha(X) \leq x$, then $h_\alpha(x) \leq 0$.
\end{proposition}

{\sc Proof.} (a) Using \eqref{EqPosNegDiff} one gets 
\[
h_\alpha(x) = \alpha \E (X - x)\One_{X > x} + (1-\alpha)\E (X - x)\One_{X \leq x}
\]
and the formula for $\psi_\alpha(x) - x$ by doing a little algebra based on the linearity of the expected value.

(b) The function $x \mapsto \E(X - x)_+$ is decreasing and the function $x \mapsto \E(X - x)_-$ is increasing. Since $0 \leq \alpha, 1 - \alpha$ the result follows.

(c) This is a consequence of (b) and $h_\alpha(e_\alpha(X)) = 0$.
\pend

\begin{corollary}
\label{CorPsi} 
One has

(a) $\psi_\alpha(x) \leq x$ if $x \geq e_\alpha(X)$;

(b) $\psi_\alpha(x) \geq x$ if $x \leq e_\alpha(X)$.
\end{corollary}

{\sc Proof.}  By \eqref{EqTwoSidedFP-Mod}, the sign of $\psi_\alpha(x) - x$ is the same as the one of $h_\alpha(x)$ since $\gamma_\alpha(x) > 0$. Thus, from Proposition \ref{PropPsiRep} (c) one gets
\[
x \geq e_\alpha(X) \; \Rightarrow \; h_\alpha(x) \leq 0 \; \Rightarrow \; \psi_\alpha(x) \leq x
\]
and
\[
x \leq e_\alpha(X) \; \Rightarrow \; h_\alpha(x) \geq 0 \; \Rightarrow \; \psi_\alpha(x) \geq x.
\]
\pend

\begin{theorem}
\label{ThmTwoSidedDescent}
Let $0 < \alpha < 1/2$. Then one has $e_\alpha(X) \leq \psi_\alpha(x)$ for each $x \in \R$.
\end{theorem}

{\sc Proof.} Abbreviate $e_\alpha = e_\alpha(X)$. First, assume $e_\alpha \leq x$. Corollary \ref{CorPsi} (a) ensures $\psi_\alpha(x) \leq x$. On the other hand,
\begin{align*}
\psi_\alpha(x) - e_\alpha 
	& = \frac{\alpha \E (X - e_\alpha)\One_{X > x} + (1-\alpha)\E (X - e_\alpha)\One_{X \leq x}}{\gamma_\alpha(x)} \\
	& = \frac{\alpha \E (X - e_\alpha)\One_{X > x} + (1-\alpha)\E (X - e_\alpha)(1 - \One_{X > x})}{\gamma_\alpha(x)} \\
	& = \frac{(2\alpha - 1)\E (X - e_\alpha)\One_{X > x} + (1-\alpha)(\E X - e_\alpha)}{\gamma_\alpha(x)} \\
	& \geq \frac{(2\alpha - 1)\E (X - e_\alpha)\One_{X > e_\alpha} + (1-\alpha)(\E X - e_\alpha)}{\gamma_\alpha(x)}.
\end{align*}
The numerator of the last fraction is
\[
(2\alpha - 1)\E (X - e_\alpha)\One_{X > e_\alpha} + (1-\alpha)(\E X - e_\alpha) = 
	\alpha \E (X - e_\alpha)\One_{X > e_\alpha} + (1-\alpha)\E (X - e_\alpha)\One_{X \leq e_\alpha} = 0
\]
since $e_\alpha$ satisfies \eqref{EqTwoSidedFP}. Altogether, $e_\alpha(X) \leq \psi_\alpha(x) \leq x$.

Secondly, assume $x \leq e_\alpha$. A similar transformation as above produces
\begin{align*}
\psi_\alpha(x) - e_\alpha 
	& = \frac{\alpha (\E X - e_\alpha) + (1 - 2\alpha) \E (X - e_\alpha)\One_{X \leq x}}{\gamma_\alpha(x)} \\
	& \geq \frac{\alpha (\E X - e_\alpha) + (1 - 2\alpha) \E (X - e_\alpha)\One_{X \leq e_\alpha}}{\gamma_\alpha(x)} = 0
\end{align*}
where the inequality is from $1 - 2\alpha \geq 0$ and $\E (X - e_\alpha)\One_{X \leq x} \geq \E (X - e_\alpha)\One_{X \leq e_\alpha}$ because of the assumptions. Thus, $e_\alpha \leq \psi_\alpha(x)$ also in this case. \pend

\begin{corollary}
\label{CorTowSidedFPConvergence}
Assume that $\{x_i\}_{i=0, 1, \ldots}$ is a sequence generated by
\begin{equation}
\label{EqPsiAlphaFixPointIteration}
x_{i+1} = \psi_\alpha(x_i), \; i = 0, 1, \ldots, 
\end{equation}
with $0 < \alpha < 1/2$ and $x_0 \in \R$ arbitrary. Then $\{x_i\}_{i=0, 1, \ldots}$ is decreasing and bounded from below by $e_\alpha(X)$ with the possible exception of $x_0$.
\end{corollary}

{\sc Proof.} If $e_\alpha \leq  x_0$, then $e_\alpha  \leq x_1 = \psi_\alpha(x_0) \leq x_0$ by Theorem \ref{ThmTwoSidedDescent}. If $e_\alpha >  x_0$, then, again by Theorem \ref{ThmTwoSidedDescent}, $e_\alpha  \leq x_1 = \psi_\alpha(x_0)$. In both cases, the remaining sequence is decreasing and bounded below by $e_\alpha$. \pend

Consequently, the sequence $\{x_i\}_{i=0, 1, \ldots}$ converges to some $\bar x$ with $e_\alpha(X) \leq \bar x \leq \max\{x_0, x_1\}$. It remains to verify $\bar x = \psi_\alpha(\bar x)  = \lim_{i \to \infty} \psi_\alpha(x_i)$.

\begin{theorem}
\label{ThmPsiContinuity} 
Under the assumption of Corollary \ref{CorTowSidedFPConvergence}, one has
\[
\bar x = \psi_\alpha(\bar x)
\]
where $\bar x = \lim_{i \to \infty} x_i$. Moreover, $\bar x = e_\alpha(X)$.
\end{theorem}

{\sc Proof.} We already know $e_\alpha \leq \psi_\alpha(x_i) \leq x_i$ for $i= 1, 2, \ldots$

First, the random variable $Y_i = (X - x_i)_+$ converges pointwise to $(X - \bar x)_+$ and is bounded since
\[
\abs{Y_i} \leq \abs{X} + \abs{x_1} =: Y \quad\text{with}\quad \E Y < +\infty.
\]
The Dominated Convergence Theorem (DCT) yields
\[
\lim_{i \to \infty} \E Y_i = \lim_{i \to \infty} \E(X - x_i)_+ = \E(X - \bar x)_+.
\]
A  parallel argument produces
\[
 \lim_{i \to \infty} \E(X - x_i)_- = \E(X - \bar x)_-.
\]

Similarly, the bounded random variables $\One_{X>x_i}$, $\One_{X \leq x_i}$ converge pointwise to $\One_{X>\bar x}$, $\One_{X \leq \bar x}$, respectively. Again by the DCT
\[
\lim_{i \to \infty} \E\One_{X>x_i} = \E\One_{X>\bar x}, \; \lim_{i \to \infty} \E\One_{X \leq x_i} = \E\One_{X \leq \bar x}.
\]
This means that the numerator as well as the denominator of $\psi_\alpha(x_i) - x_i$ in \eqref{EqTwoSidedFP-Mod} converge to the corresponding expressions for $\bar x$. Moreover, the denominator can only have values in $[\alpha, 1-\alpha]$ and is thus bounded away from 0. This gives
\[
\lim_{i \to \infty}\psi_\alpha(x_i) = \psi_\alpha(\bar x).
\]
Since both sides of \eqref{EqPsiAlphaFixPointIteration} converge to the same limit $\bar x$, the claim follows: the right hand side of \eqref{EqTwoSidedFP-Mod} yields $\psi_\alpha(\bar x) - \bar x = 0$ and thus $\bar x = e_\alpha$ since $e_\alpha$ is the only root of the function $\psi_\alpha(x) - x$. \pend

\begin{remark}
Of course, Theorems \ref{ThmTwoSidedDescent} and \ref{ThmPsiContinuity} have counterparts for the case $1/2 < \alpha < 1$. However, one can also use the formula $e_\alpha(X) = -e_{1-\alpha}(-X)$ to transform the problem. 
\end{remark}

\subsection{Fixed point iterations for sample expectiles}

Assume that a sample of data observations $Z = \{z_1, z_2, ..., z_N\} \subseteq \R$ is available. Moreover, it is assumed that they are numbered in increasing order, i.e., $z_1 \leq \ldots \leq z_N$. The empirical sample expectile $e_{N, \alpha}$ is the fixed point of 
\begin{equation}
\label{EqEmpiricalExpectile}
\hat{\varphi}_\alpha(x) =  
\left \{ 
\begin{aligned} 
& \bar z + \frac{2\alpha - 1}{N(1-\alpha)} \sum_{z_n \geq x} (z_n - x) \quad & \text{for} \quad 0 < \alpha < \frac{1}{2} \\
& \bar z + \frac{2\alpha - 1}{N \alpha} \sum_{z_n < x} (x - z_n ) \quad & \text{for} \quad \frac{1}{2} < \alpha < 1 
\end{aligned} 
\right.
\end{equation}
where $\bar z = 1/N \sum_{n=1}^N z_n$ is the sample mean and $e_{N, 0.5} = \bar z$ (see also, for example, \cite[formula (8)]{CascosOchoa21JMVA}). The fixed point iteration for $\hat{\varphi}_\alpha$ terminates after a finite number of steps. We start with a preliminary result.

\begin{proposition}
If $x \in [z_1, z_N]$, then $\hat{\varphi}_\alpha(x) \in [z_1, z_N]$.
\end{proposition}

{\sc Proof.} We consider the case $0 < \alpha < \frac{1}{2}$. On the one hand, one has $2\alpha - 1 < 0$ which implies $\hat{\varphi}_\alpha(x) \leq \bar z \leq z_N$. On the other hand, one has
\[
\sum_{z_n \geq x} (z_n - x) \leq \sum_{z_n \geq x} (z_n - z_1) \leq \sum_{n=1}^N (z_n - z_1) = N\bar z - Nz_1.
\]
This implies
\[
\hat{\varphi}_\alpha(x) \geq \bar z +\frac{2\alpha - 1}{N(1-\alpha)}N(\bar z - z_1) = \bar z + \frac{2\alpha - 1}{1-\alpha}(\bar z - z_1) =
		\frac{\alpha}{1-\alpha}\bar z + \frac{1 - 2\alpha}{1-\alpha}z_1
\]
which means that $\hat{\varphi}_\alpha(x)$ is greater or equal than a convex combination of $\bar z$ and $z_1$, hence $\hat{\varphi}_\alpha(x) \geq z_1$. Altogether, $\hat{\varphi}_\alpha(x) \in [z_1, \bar z] \subseteq [z_1, z_N]$. 

The proof for the case $\frac{1}{2} < \alpha < 1 $ follows mirrored lines and yields $\hat{\varphi}_\alpha(x) \in [\bar z, z_N] \subseteq [z_1, z_N]$.\pend

Note that $\hat{\varphi}_\alpha(x) \leq \bar z$ for $0 < \alpha < \frac{1}{2}$ and $\bar z \leq \hat{\varphi}_\alpha(x)$ for $\frac{1}{2}< \alpha < 1$ also is a consequence of the monotonicity of the $\alpha$-expectile with respect $\alpha$.

\begin{theorem}
Let $0 < \alpha < \frac{1}{2}$ and $x_0 \in [z_1, z_N]$. The fixed point iteration
\begin{equation}
\label{EqFiniteFPI}
x_{i+1} = \hat{\varphi}_\alpha(x_i), \; i = 0, 1, \ldots 
\end{equation}
terminates after a finite number of steps in the following sense: if $z_K < e_{N, \alpha} < z_{K+1}$ for $K \in \{1, \ldots, N-1\}$, then
there is $i \in \N\bs\{0\}$ such that
\begin{equation}
\label{EqStablePartition}
\forall j = 0, 1, \ldots \colon \{z \in Z \mid z > x_{i + j}\} = \{z \in Z \mid z > e_{N, \alpha}\} = \{z_{K+1}, \ldots, z_N\}
\end{equation}
in which case
\begin{equation}
\label{EqFiniteTerminationFormula}
e_{N, \alpha} = \frac{(2\alpha - 1)\tilde z + N(1 - \alpha)\bar z}{N(1 - \alpha)+(N-K)(2\alpha - 1)} 
\end{equation}
where $\tilde z = \sum_{k = K+1}^N z_k$. 
\end{theorem}

{\sc Proof.} The Banach Fixed Point Theorem \cite[3.48 Contraction Mapping Theorem]{AliprantisBorder06Book} provides the estimate
\[
\forall i = 0, 1, \ldots \colon \norm{x_i - e_{N, \alpha}}_E \leq \gamma^i \norm{x_0 - e_{N, \alpha}}_E
\]
where $\gamma = \frac{1 - 2\alpha}{\alpha}$ is the contraction constant and $\norm{\cdot}_E$ the euclidean norm. Thus, for $i$ large enough, $\norm{x_i - e_{N, \alpha}}_E < \min\{\norm{z_K - e_{N, \alpha}}_E, \norm{z_{K+1} - e_{N, \alpha}}_E\}$, and this inequality remains true for $x_{i+j}$, $j = 1,2, \ldots$. This implies \eqref{EqStablePartition}. Now, one can solve the linear equation
\[
e_{N, \alpha} = \bar z + \frac{2\alpha - 1}{N(1-\alpha)} \sum_{n = K+1}^N (z_n - e_{N, \alpha}) 
\]
for $e_{N, \alpha}$ which yields \eqref{EqFiniteTerminationFormula}. \pend

One may wonder what happens if $e_{N, \alpha}$ coincides with a data point, i.e.,  $e_{N, \alpha} \in Z$. Here is a simple example.
\begin{example}
\label{ExThreePoint} 
Let us consider the data set $Z = \{1,2,7\}$ and $\alpha = 1/6$. One has $\bar z = 10/3$ and $e_{3, \alpha} = 2$ since $\hat{\varphi}_\alpha(2) = 2$ as one may see from
\[
\hat{\varphi}_\alpha(x) = 
	\left\{
	\begin{array}{ccc}
	\frac{10}{3} & \colon & x \geq 7 \\[.1cm]
	\frac{22}{15} + \frac{4}{15}x & \colon & 2 \leq x < 7 \\[.1cm]
	\frac{14}{15} + \frac{8}{15}x & \colon & 1 \leq x < 2  \\[.1cm]
	\frac{2}{3} +  \frac{4}{5}x & \colon & x < 1
	\end{array}
	\right.
\]
Starting the fixed point iteration for $\hat{\varphi}_\alpha$ with $x_0 = \bar z = 10/3$ one gets $x_1 = 106/45 \approx 2.355556$, $x_2 = 1414/675 \approx 2.094815$, $x_3 = 20506/10125 \approx 2.025284$. These values are decreasing and greater than $e_{3, \alpha} = 2$. Starting with $x_0 = 1 = \min Z$ one gets $x_1 = 22/15 \approx 1.466667$, $x_2 = 386/225 \approx 1.715556$, $x_3 = 6238/3375 \approx 1.848296$. These values are increasing and less than $e_{3, \alpha} = 2$. In both cases, the set $\cb{z \in Z \mid z > x_i}$ does not change anymore starting with $x_1$. Thus, one can compute $e_{3, \alpha}$ via \eqref{EqStablePartition} as
\[
e_{3, 1/6} = \frac{(-2/3)9 + 3 \cdot 5/6 \cdot 10/3}{3 \cdot 5/6 +(3-1) \cdot (-2/3)} = 2.
\] 
\end{example}

The following result shows that the behaviour in Example \ref{ExThreePoint} is not a coincidence. 

\begin{proposition}
\label{PropPhiHeadFeatures} Let $0 < \alpha < \frac{1}{2}$.

(a) The function $\hat{\varphi}_\alpha$ is piecewise linear and increasing with $\hat{\varphi}_\alpha(x) = \bar z$ for all $x \geq z_N$.

(b) If $x \leq e_{N, \alpha}$ then $\hat{\varphi}_\alpha(x) \leq e_{N, \alpha}$; if $x \geq e_{N, \alpha}$ then $\hat{\varphi}_\alpha(x) \geq e_{N, \alpha}$.
\end{proposition}

{\sc Proof.} (a) The monotonicity follows from the monotonicity of the positive part function and $\frac{2\alpha - 1}{N(1-\alpha)} < 0$. Piecewise linearity is a consequence of the definition of $\hat{\varphi}_\alpha$. (b) This follows from (a) and $\hat{\varphi}_\alpha(e_{N, \alpha}) = e_{N, \alpha}$. \pend

The previous proposition implies that the sequence defined by \eqref{EqFiniteFPI} stays above and below $e_{N, \alpha}$ if the starting point $x_0$ is above and below $e_{N, \alpha}$, respectively. Moreover, the sequence is  decreasing and increasing, respectively.

By the way of conclusion, \eqref{EqStablePartition} remains valid with $z_{K+1} = e_{N, \alpha}$. 

\medskip
Using the function $\psi_\alpha$ instead of $\varphi_\alpha$, the empirical sample expectile is the fixed point of the function
\begin{equation}
\label{EqEmpiricalExpectileFP2}
\hat\psi_\alpha(x) = \frac{\alpha\sum_{z \in Z, z > x} z + (1-\alpha)\sum_{z \in Z, z \leq x}z}
{\alpha \sum_{\{z \in Z \mid z>x\}} \mathbbm{1}_{\{z\}}+ (1-\alpha)\sum_{\{z \in Z \mid z \leq x\}} \mathbbm{1}_{\{z\}}}.
\end{equation}

Theorem 2.2 in \cite{DaouiaStupflerUsseglio24SC} provides a result on quadratic convergence for the fixed point iteration for $\hat\psi_\alpha$ for continuous distributions, not for a discrete data set. In the discrete case, the convergence follows from the general case in Theorem \ref{ThmPsiContinuity}. 

A few features of $\hat\psi_\alpha$ simplify the computations.

\begin{theorem}
\label{ThmDiscreteTwoSidedFeatures}
The function $\hat\psi_\alpha$ is piece-wise constant. Namely, if $Z = \{z_1, \ldots, z_N\}$ with $z_1 < \ldots < z_N$, then
\begin{equation}
\label{EqEmpiricalExpectileFP3}
\forall x \in (z_n, z_{n+1}] \colon 
	\hat\psi_\alpha(x) = \frac{\alpha\sum_{k = n + 1}^N z_n + (1-\alpha)\sum_{k=1}^n z_n}{\alpha n + (1-\alpha)(N-n)}.
\end{equation}
The function $x \mapsto \hat\psi_\alpha(x) - x$ is piece-wise linear and piece-wise strictly decreasing. It can only jump up at points $z_n \in Z$ with $e_{N, \alpha} \leq z_n$.
\end{theorem}

{\sc Proof.} The first part is an immediate consequence of the definition of $\hat\psi_\alpha$ while the second follows from the first and Proposition \ref{PropPsiRep} (b), (c). \pend

\begin{figure}[ht]
\includegraphics[width=0.7\textwidth]{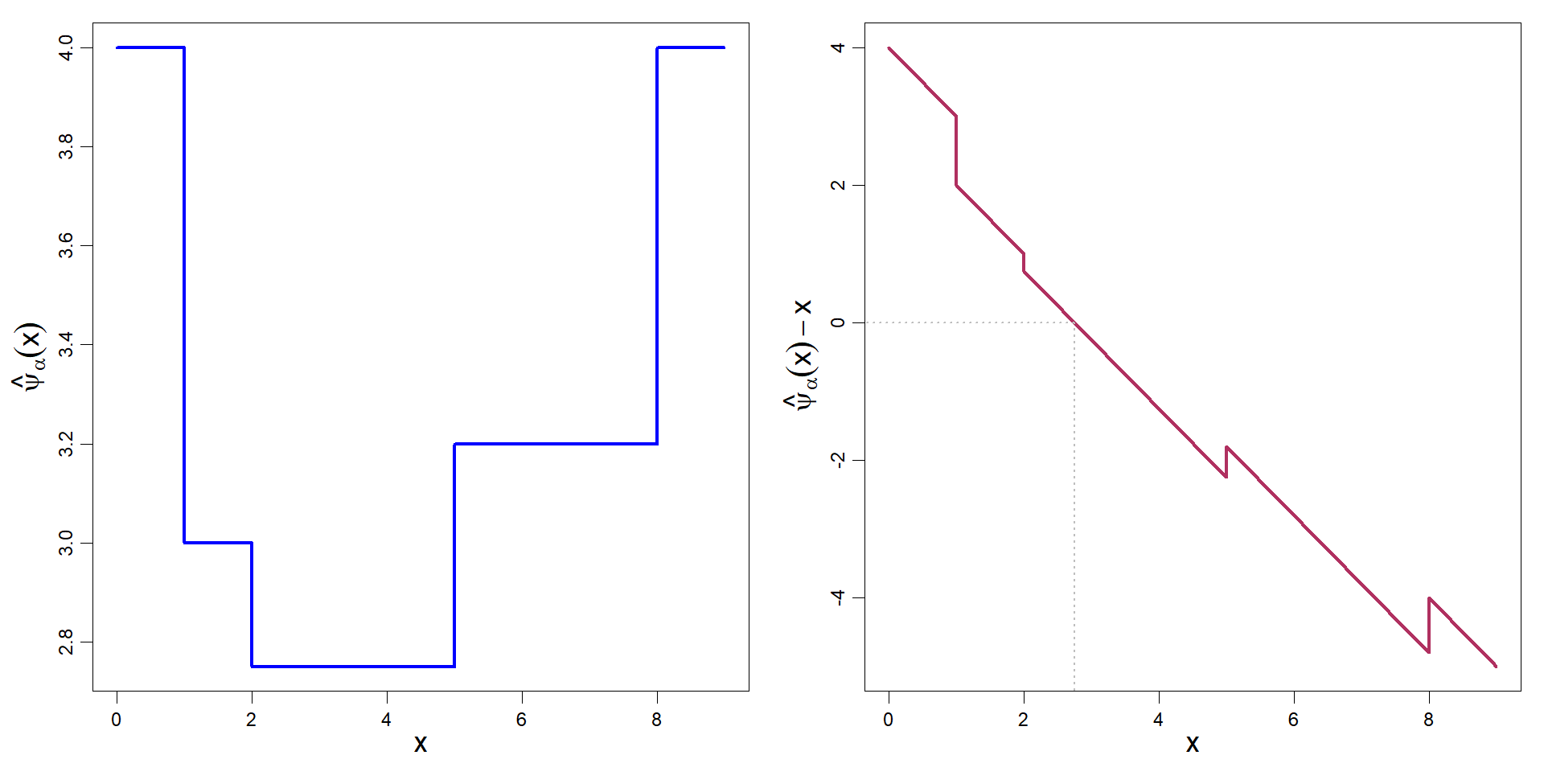}
\caption{Function $\hat{\psi}_\alpha(x) - x$ from Example \ref{Eg:Jump}.}
\label{Fig_Jump}
\end{figure}

\begin{example}\label{Eg:Jump}
Consider $Z = \{1, 2, 5, 8\}$ and $\alpha = 1/4$. One has $\bar z = 4$ and $e_{4, 0.25} = 11/4 = 2.75$ since $\hat{\psi}_\alpha(2.75) = 2.75$ as one may see from
\[
\hat{\psi}_\alpha(x) = 
	\left\{
	\begin{array}{ccc}
	4 & \colon & x < 1 \\[.1cm]
	3 & \colon & 1 \leq x < 2 \\[.1cm]
	\frac{11}{4} & \colon & 2 \leq x < 5  \\[.1cm]
	\frac{16}{5}  & \colon & 5 \leq x < 8 \\[.1cm]
	4 & \colon & 8 \leq x  
	\end{array}
	\right.
\]
The function $x \mapsto \hat{\psi}_\alpha(x)$ is piece-wise constant and $\hat{\psi}_\alpha(x) - x$ "jumps up" at $x = 5$ and again at $x = 8$ since $h_\alpha(5), h_\alpha(5) < 0$: see Figure \ref{Fig_Jump} (right).
\end{example}

\begin{corollary}
Assume $e_{N, \alpha} < x \in (z_n, z_{n+1}]$. Then
\[
\hat\psi_\alpha(x) \leq z_n \quad \text{or} \quad \hat\psi_\alpha(x) = e_{N, \alpha}.
\]
\end{corollary}

{\sc Proof.}  One has $\hat\psi_\alpha(x) \leq x$ by Corollary \ref{CorPsi} (a). Assume $\hat\psi_\alpha(x) > z_n$. Then $\hat\psi_\alpha(x) = \hat\psi_\alpha(\hat\psi_\alpha(x))$ since $\hat\psi_\alpha$ is constant on $(z_n, z_{n+1}]$. Hence $\hat\psi_\alpha(x) = e_{N, \alpha}$. Otherwise, $\hat\psi_\alpha(x) \leq z_n$. \pend

{\bf Finite termination.} By the way of conclusion, the sequence $x_{i+1} = \hat\psi_\alpha(x_i)$ either stays in an interval $(z_n, z_{n+1}]$ (i.e., $x_i, x_{i+1} \in (z_n, z_{n+1}]$) in which case $x_{i+1} = \hat\psi_\alpha(x_i) = e_{N, \alpha}$ or it "moves down" at least one interval: from $(z_n, z_{n+1}]$ to $(z_{n-\ell}, z_{n - \ell +1}]$ with $\ell \geq 1$. Thus, the fixed point iteration for $\hat\psi_\alpha$ terminates after at most $N$ steps. 

\medskip
The following example illustrates this behavior.

\begin{figure}
\includegraphics[width=0.7\textwidth]{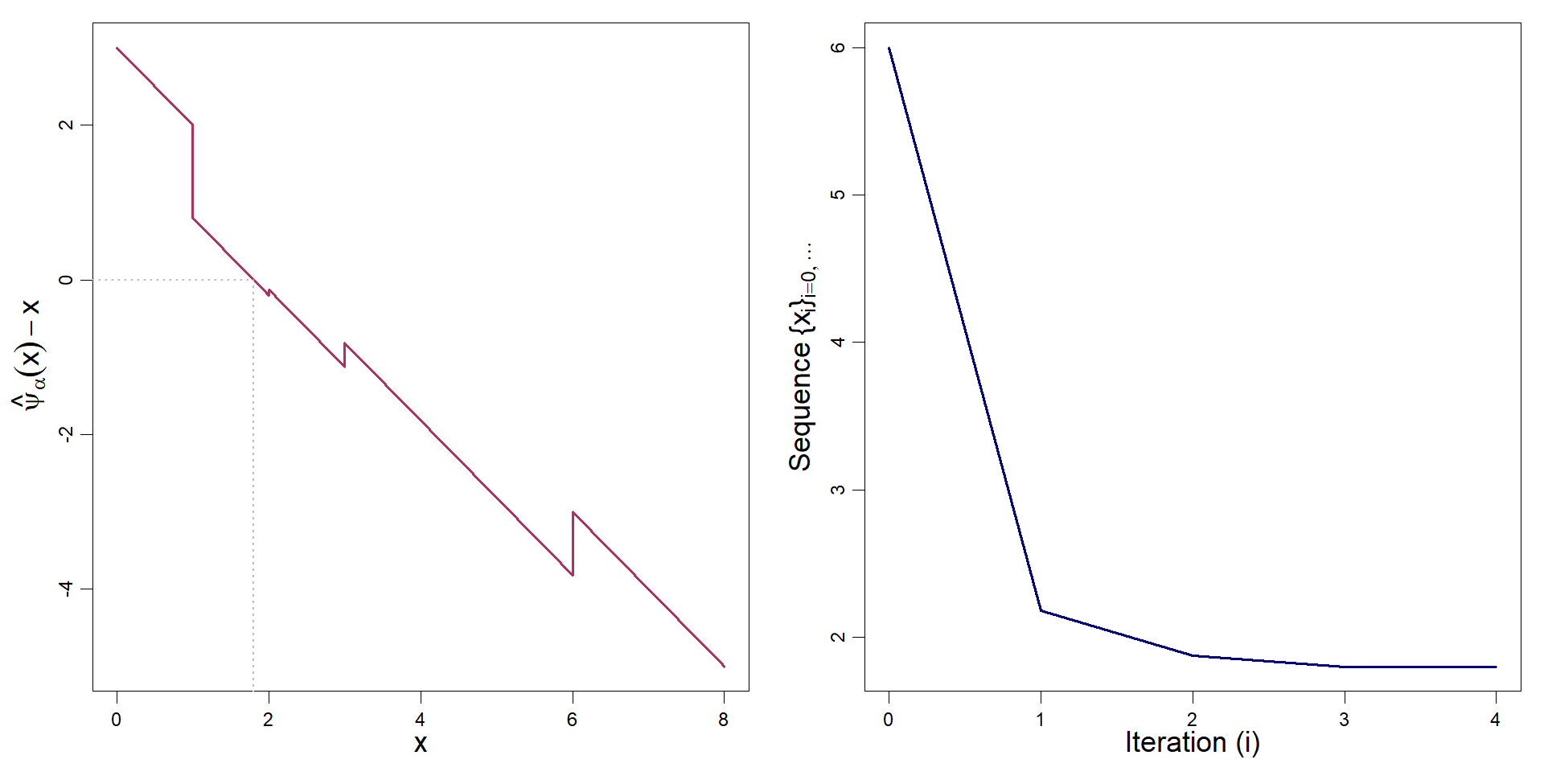}
\caption{Illustration of Example \ref{Eg:TerminalIT}}
\label{Fig_TerminalIT}
\end{figure}

\begin{example}
\label{Eg:TerminalIT}
Consider $Z = \{1, 2, 3, 6\}$ and $\alpha = 1/8$. One has $\bar z = 3$ and $e_{4, 0.125} = 1.8$ since $\hat{\psi}_\alpha(1.8) = 1.8$ as one may see from
\[
\hat{\psi}_\alpha(x) = 
	\left\{
	\begin{array}{ccc}
	3 & \colon & x < 1 \\[.1cm]
	\frac{9}{5} & \colon & 1 \leq x < 2 \\[.1cm]
	\frac{15}{8} & \colon & 2 \leq x < 3  \\[.1cm]
	\frac{24}{11}  & \colon & 3 \leq x < 6 \\[.1cm]
	3 & \colon & 6 \leq x  
	\end{array}
	\right.
\]
Again, the function $x \mapsto \hat{\psi}_\alpha(x)$ is piece-wise constant and $\hat{\psi}_\alpha(x) - x$ strictly decreasing except at $x = 2$, $x = 3$, $x = 6$. The fixed point iteration for $\hat{\psi}_\alpha$ starting with $x_0 = 6$ produces $x_1 = \bar z = 3$, $x_2 = 24/11$, $x_3 = 15/8$ and $x_4 = 9/5 = 1.8 = e_{4, 0.125}$. Thus, in this case, $N = 4$ steps are needed to compute the expectile: see Figure \ref{Fig_TerminalIT}.
\end{example}

Both of these behaviors are empirically illustrated in Figure \ref{AlgoCompare} below. The data set for the underlying example is from a random sample of 1000 data points drawn from the normal distribution with mean $\mu = 2$ and standard deviation $\sigma = 6$. One may see that even for such a large data set, the convergence/finite termination can happens rather fast. 

The two leftmost diagram rows also illustrate the "overshooting effect" in the first step according to Corollary \ref{CorTowSidedFPConvergence}. For the two middle pictures, the $\alpha$-quantile with the same $\alpha$ as for the expectile was used as starting point which illustrates that "overshooting" may also occur for $\min Z < x_0$. 

The impression that the pink curve seems to be non-concave (upper row, left and middle picture) is due to the fact that the finite termination criterion \eqref{EqFiniteTerminationFormula} cuts off the true graph of $\hat \varphi_\alpha$ which indeed is a concave function (compare formula \eqref{EqEmpiricalExpectile}).

\begin{figure}
\includegraphics[width=0.9\textwidth]{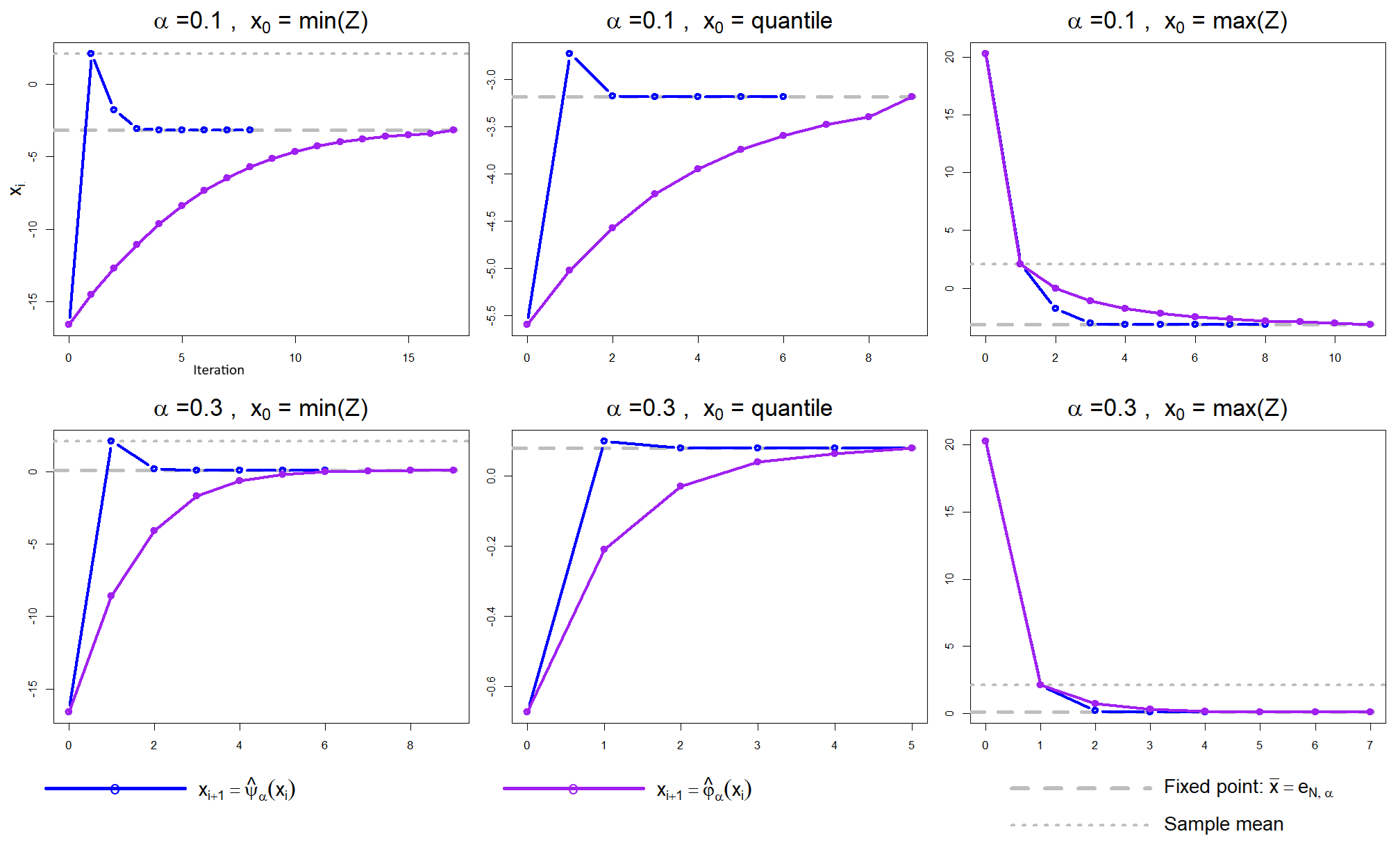}
\caption{Convergence of two formulas to the fixed point}
\label{AlgoCompare}
\end{figure}

\section{Conclusion}

Expectiles are minimizers of asymmetric quadratic loss functions. The corresponding first order condition can be turned into fixed point form by either replacing the positive or the negative part by the respective other one or, alternatively, by keeping both parts. This leads to different fixed point characterizations of expectiles, a one-sided one \eqref{EQOneSidedFP-A}, \eqref{EQOneSidedFP-B} and a two-sided one \eqref{EqTwoSidedFP}. Existence and uniqueness of the fixed point can be shown via Banach's fixed point theorem applied to the one-sided version. The question remains if such a procedure can be designed also for other statistical parameters, especially M-quantiles \cite{BelliniEtAl14IME}.

Both fixed point versions can be exploited for computational procedures which was already implemented in \cite{SobotkaEtAl14Github}, but a general convergence proof was not available yet. This paper closes this little gap in the literature providing also finite termination criteria for both the one-sided and the two-sided fixed point iteration for sample expectiles.



\bibliographystyle{plain}

\end{document}